\documentclass[a4paper,11pt]{article}

\usepackage{amsmath,amsthm,amsfonts,amssymb,amscd}
\usepackage{verbatim}
\usepackage{color}
\usepackage{hyperref,cite,mdwlist}
\usepackage{manyfoot}
\usepackage[english]{babel}
\usepackage{cite}
\usepackage{latexsym}
\usepackage{amsopn}
\usepackage[all]{xy}
\usepackage{mathrsfs}
\usepackage{float}
\usepackage[active]{srcltx}
\usepackage{fullpage}
\usepackage{tikz}
\usepackage{dsfont}
\usetikzlibrary{matrix}

\newtheorem{inizio}{Lemma}[section]
\newtheorem{theorem}[inizio]{Theorem}

\newtheorem*{theorem*}{Theorem}

\theoremstyle{definition}
\newtheorem{definition}[inizio]{Definition}
\newtheorem{example}[inizio]{Example}

\renewcommand{\to}{\longrightarrow}
\newcommand{\lr}{\longrightarrow}

\title{Isotrivially fibred surfaces and their numerical invariants}

\author{Francesco Polizzi}

\date{}

\begin{document}

\maketitle

\begin{abstract}
We give a survey of our previous work on relatively minimal
isotrivial fibrations $\alpha \colon X \to C$, where $X$ is a
smooth, projective surface and $C$ is a curve. In particular, we
consider two inequalities involving the numerical invariants
$K_X^2$ and $\chi(\mathcal{O}_X)$ and we illustrate them by means
of several examples and counterexamples.
\end{abstract}

\Footnotetext{{}}{\textit{2010 Mathematics Subject
Classification}: 14J99, 14J29}

\Footnotetext{{}} {\textit{Keywords}: isotrivial fibrations,
numerical invariants}







\section{Introduction} \label{sec:intro}

One of the most useful tools in the study of algebraic surfaces is
the analysis of \emph{fibrations}, that is morphisms with
connected fibres from a smooth surface $X$ to a curve $C$. A
fibration $\alpha \colon X \to C$ is called \emph{isotrivial} when
all its smooth fibres are isomorphic; a deep investigation of such
kind of fibrations can be found in \cite{Ca00}.

 A special kind of isotrivial fibrations are
the standard ones, whose definition is as follows. A smooth
surface $X$ is called a \emph{standard isotrivial fibration} if
there exists a finite group $G$, acting faithfully on two smooth
curves $C_1$ and $C_2$ and diagonally on their product, such that
$X$ is isomorphic to the minimal desingularization of $T=(C_1
\times C_2)/G$. Then $X$ has two isotrivial fibrations $\alpha_i
\colon X \to C_i/G$, induced by the natural projections $p_i
\colon C_1 \times C_2 \to C_i$.

If the action of $G$ on $C_1 \times C_2$ is free, we say that $X=T$
is a \emph{quasi-bundle}. For instance, unmixed Beauville surfaces
are defined as rigid quasi-bundles, i.e., quasi-bundles such that
$C_i/G \cong \mathbb{P}^1$ and the covers $C_i \to \mathbb{P}^1$ are
both branched at three points. Therefore the classification of
Beauville surfaces can be reduced to combinatorial finite group
theory involving triangle groups; we refer the reader to the other
papers in this Volume for further details.

Standard isotrivial fibrations were thoroughly investigated by F.
Serrano in \cite{Se90} and \cite{Se96}; in particular he showed, by
a monodromy argument, that every isotrivial fibration is isomorphic
to a standard one. Since then, such fibrations have been widely used
in order to produce new examples of minimal surfaces of general type
with small birational invariants, in particular with $p_g=q=0$
(\cite{BaCa04}, \cite{BaCaGr08}, \cite{BaCaPi11}, \cite{BaCaGrPi12},
\cite{BaPi12}), with $p_g=q=1$ (\cite{Pol06}, \cite{CarPol07},
\cite{Pol08}, \cite{Pol09}, \cite{MiPol10}, \cite{Pol10} ) and with
$p_g=q=2$ (\cite{Pe11}).

In this paper we discuss the following theorem, that was  obtained
in \cite{Pol10}. Let $f \colon X \to C$ be any relatively minimal,
isotrivial fibration with $g(C) \geq 1$. If $X$ is neither ruled nor
isomorphic to a quasi-bundle, then $K_X^2 \leq 8 \chi
(\mathcal{O}_X)-2$. If, in addition, $K_X$ is ample, then $K_X^2
\leq 8 \chi (\mathcal{O}_X)-5$. This generalizes previous results of
Serrano and Tan (\cite{Se96}, \cite{Tan96}).

This work is organized as follows. In Section \ref{sec:main-result}
we set up the notation and the terminology, we state our theorem and
we provide a sketch of its proof. In Section \ref{sec:examples} we
exhibit several examples and counterexamples illustrating its
meaning. More precisely, Examples \ref{ex.1} and \ref{ex.2} imply
that both the above inequalities involving $K_X^2$ and
$\chi(\mathcal{O}_X)$ are sharp, whereas Examples \ref{ex.3} and
\ref{ex.4} show that, when $K_X$ is not ample, both cases $K_X^2 = 8
\chi (\mathcal{O}_X)-3$ and $K_X^2 = 8 \chi (\mathcal{O}_X)-4$
actually occur. Finally, in Example \ref{ex.5} we describe an
isotrivially fibred surface $X$ with $K_X^2 = 8 \chi
(\mathcal{O}_X)-5$ and $K_X$ \emph{not} ample. In all these examples
$X$ is a minimal surface of general type with $p_g=q=1$, obtained as
a standard isotrivial fibration.

\bigskip
$\mathbf{Notations \; and \; conventions}$. All varieties in this
article are defined over $\mathbb{C}$. If $X$ is a projective,
non-singular surface $X$ then $K_X$ denotes the canonical class,
$p_g(X)=h^0(X, \, K_X)$ is the \emph{geometric genus}, $q(S)=h^1(X,
\, K_X)$ is the \emph{irregularity} and
$\chi(\mathcal{O}_X)=1-q(X)+p_g(X)$ is the \emph{Euler
characteristic}.

If $G$ is a finite group and $g \in G$, we denote by $|G|$ and
$o(g)$ the orders of $G$ and $g$, respectively.

\bigskip
$\mathbf{Acknowledgements.}$ This paper is an expanded version of
the talk given by the author at the conference \emph{Beauville
surfaces and Groups}, Newcastle University (UK), 7-9th June 2012.
The author is grateful to the organizers N. Barker, I. Bauer, S.
Garion and A. Vdovina for the invitation and the kind hospitality.
He was partially supported by Progetto MIUR di Rilevante Interesse
Nazionale \emph{Geometria delle Variet$\grave{a}$ Algebriche e loro
Spazi di Moduli}. He also thanks the referee, whose comments helped
to improve the presentation of these results.

\section{The main result} \label{sec:main-result}

\begin{definition} \label{def:iso-fib}
Let $X$ be a smooth, complex projective surface and let $\alpha
\colon X \to C$ be a fibration onto a smooth curve $C$. We say
that $\alpha$ is \emph{isotrivial} if all its smooth fibres are
isomorphic.
\end{definition}

\begin{definition} \label{def:stand-iso}
A smooth surface $S$ is called a \emph{standard isotrivial
fibration} if there exists a finite group $G$, acting faithfully on
two smooth projective curves $C_1$ and $C_2$ and diagonally on their
product, so that $S$ is isomorphic to the minimal desingularization
of $T:=(C_1 \times C_2)/G$. We denote such a desingularization by
$\lambda \colon S \to T$.
\end{definition}
If $\lambda \colon S \to T=(C_1 \times C_2)/G$ is any standard
isotrivial fibration, composing the two projections $\pi_1 \colon T
\lr C_1/G$ and $\pi_2 \colon T \lr C_2/G$ with $\lambda$ one obtains
two morphisms $\alpha_1 \colon S \lr C_1/G$ and $\alpha_2 \colon S
\lr C_2/G$, whose smooth fibres are isomorphic to $C_2$ and $C_1$,
respectively. One also has $q(S)=g(C_1/G) +g(C_2/G)$, see
\cite{Fre71}.

A monodromy argument (\cite[Sect. 2]{Se96}) implies that any
isotrivial fibration $\alpha \colon X \to C$ is birational to a
standard one; in other words, there exists $T= (C_1  \times C_2)/G$
and a birational map $T \dashrightarrow X$ such that the following
diagram
\begin{equation} \label{diagram-iso}
\xymatrix{ S \ar[d]_{\lambda}  \ar[rd]^{\psi} \\
T \ar@{-->}[r]  \ar[d]_{\pi_2} &
 X \ar[d]^{\alpha} \\
 C_2/G \ar[r]^{\cong} & C}
\end{equation}
commutes.

When the action of $G$ is free, then $S=T$ is called a
\emph{quasi-bundle}; in this case one has $K_S^2 = 8 \chi
(\mathcal{O}_S)$. In 1996, F. Serrano and, indipendently, S. L. Tan
improved this result, showing that for any isotrivial fibration
$\alpha \colon X \to C$ one has
\begin{equation} \label{eq:serrano-tan}
K_X^2 \leq 8 \chi (\mathcal{O}_X)
\end{equation}
and that the equality holds if and only if $X$ is either ruled or
isomorphic to a quasi bundle (\cite{Se96}, \cite{Tan96}).
Serrano's proof is based on a fine analysis of the projective
bundle $\mathbb{P}(\Omega^1_X)$, whereas Tan's proof uses base
change techniques.

This paper deals with the following refinement of
\eqref{eq:serrano-tan}, that we proved in \cite{Pol10}.

\begin{theorem} \label{thm:main}
Let $\alpha \colon X \to C$ be any relatively minimal isotrivial
fibration, with $g(C) \geq 1$. If $X$ is neither ruled nor
isomorphic to a quasi bundle, then
\begin{equation} \label{eq:main-1}
K_X^2 \leq 8 \chi (\mathcal{O}_X) -2,
\end{equation}
and if the equality holds then $X$ is a minimal surface of general
type whose canonical model has precisely two ordinary double points
as singularities. Moreover, under the further assumption that $K_X$
is ample, we have
\begin{equation} \label{eq:main-2}
K_X^2 \leq 8 \chi (\mathcal{O}_X) -5.
\end{equation}
Finally, inequalities \eqref{eq:main-1} and \eqref{eq:main-2} are
sharp.
\end{theorem}
Let us give a sketch of the proof of Theorem \ref{thm:main}, whose
full details can be found in \cite{Pol10}. Let us consider a
standard isotrivial fibration $S$ birational to $X$ and the
corresponding commutative diagram \eqref{diagram-iso}. Since
$\alpha$ is relatively minimal and $g(C) \geq 1$, the surface $X$
is a minimal model. Moreover we are assuming that $X$ is not
ruled, so $K_X$ is nef and the birational map $\psi \colon S \to
X$ is actually a morphism (\cite[Proposition 8]{Pet04}), which
induces an isomorphism of $X$ with the minimal model $S_m$ of $S$.
Since $C_2/G$ has positive genus, all the $(-1)$-curves of $S$ are
necessarily contained in fibres of $\alpha_2 \colon S \to C_2 /G$,
hence
 our isotrivial fibration $\alpha  \colon X \to C$ is equivalent to
 an isotrivial fibration $\alpha_m \colon S_m \to C_2/G$.

The surface $T$ contains at most a finite number of isolated
singularities that, locally analytically, look like the quotient of
$\mathbb{C}^2$ by the action of the cyclic group $\mathbb{Z}/n
\mathbb{Z}=\langle \xi \rangle$ defined by $\xi \cdot (x, \,y)=(\xi
x, \, \xi^qy)$, where $0 < q < n$, $(n, \, q)=1$ and
 $\xi$  is a primitive $n$-th root of unity. We call this
singularity a \emph{cyclic quotient singularity} of type
$\frac{1}{n}(1, \,q)$.  The exceptional divisor $\mathcal{E}$ of its
minimal resolution is a HJ-string
 (abbreviation of Hirzebruch-Jung string), that is to say, a
 connected union $\mathcal{E}=\bigcup_{i=1}^k Z_i$ of smooth rational
 curves $Z_1, \ldots, Z_k$ with
 self-intersection $-b_i:= Z_i^2\leq -2$, and ordered linearly so that $Z_i
 Z_{i+1}=1$ for all $i$, and $Z_iZ_j=0$ if $|i-j| \geq 2$.
More precisely, given the continued fraction
\begin{equation*}
\frac{n}{q}=[b_1,\ldots,b_k]=b_1-
                                \cfrac{1}{b_2 -\cfrac{1}{\dotsb
                                 - \cfrac{1}{\,b_k}}}, \quad b_i\geq 2 ~,
\end{equation*}
the dual graph of $\mathcal{E}$ is
{\setlength{\unitlength}{1.1cm}
\begin{center}
\begin{picture}(1,0.5)
\put(0,0){\circle*{0.2}} \put(1,0){\circle*{0.2}}
\put(0,0){\line(1,0){1}} \put(-0.3,0.2){\scriptsize $-b_1$}
\put(0.70,0.2){\scriptsize $-b_2$} \put(2,0){\circle*{0.2}}
\put(1,0){\line(1,0){0.2}} \put(1.3,0){\line(1,0){0.15}}
\put(1.55,0){\line(1,0){0.15}} \put(1.8,0){\line(1,0){0.2}}
\put(3,0){\circle*{0.2}} \put(2,0){\line(1,0){1}}
\put(1.70,0.2){\scriptsize $-b_{k-1}$} \put(2.70,0.2){\scriptsize
$-b_k$}
\end{picture}         \hspace{2.5cm}
\end{center}
} \vspace{.5cm} For instance, the cyclic quotient singularities
$\frac{1}{n}(1, \, n-1)$ are precisely the rational double points of
type $A_{n-1}$; in particular, the singularities $\frac{1}{2}(1, \,
1)$ are the ordinary double points.

Finally, the invariants $K_{S_m}^2$ and $e(S_m)$ can be computed
knowing the number and type of the singularities of $T$. In fact,
since $g(C_2/G)=g(C) \geq 1$, all the $(-1)$-curves of $S$ are
components of reducible fibres of $\alpha_2 \colon S \to C_2 /G$; it
follows that it is possible to define, for any such a reducible
fibre $F$, an invariant
 $\delta(F) \in \mathbb{Q}$ such that
 \begin{equation} \label{aaaaa}
 K_{S_m}^2 = 8 \chi (\mathcal{O}_{S_m})- \sum_{F \;
\textrm{reducible} } \delta(F).
\end{equation}
The proof of Theorem \ref{thm:main} follows from a careful
analysis of the possible values of $\delta(F)$, based on some
identities on continued fractions which are a consequence of the
so-called Riemenschneider's duality between the HJ-expansions of
$\frac{n}{q}$ and $\frac{n}{n-q}$.

\section{Examples} \label{sec:examples}

Let us denote by $\Gamma(0 \;|\; m_1, \ldots, m_r)$ the abstract
group of Fuchsian type with presentation
\begin{equation*}
\begin{split}
\Big \langle g_1, \ldots, g_r \; | \; g_k^{m_k}=1, \; \;
\prod_{i=1}^r g_i=1 \Big \rangle
\end{split}
\end{equation*}
and by $\Gamma(1 \;|\; n_1, \ldots, n_s)$ the abstract group of
Fuchsian type with presentation
\begin{equation*}
\begin{split}
\Big \langle l_1, \ldots, l_s, \; h_1, \, h_2 \; | \; l_k^{n_k}=1,
\; \; [h_1, \, h_2] \prod_{i=1}^s l_i=1 \Big \rangle.
\end{split}
\end{equation*}
We call the sets  $(m_1, \ldots, m_r)$ and $(n_1, \ldots, n_s)$ the
$\emph{branching data}$ of these groups. For convenience we make
abbreviations such as $(2^3,3^2)$ for $(2,\,2,\,2,\,3,\,3)$
 when we write down the branching data.

\begin{example} \label{ex.1} This example appears in
\cite[Section 7]{Pol09}.

Take  $G= D_8$, the dihedral group of order $8$, with presentation
\begin{equation*}
G= \langle x, \, y \, | \, x^2=y^4=1, \, xy=y^3x \rangle.
\end{equation*}
There are two epimorphisms of groups
\begin{equation*}
\varphi \colon \Gamma(0 \, | \, 2^4, \, 4) \to G, \quad \psi
\colon \Gamma(1 \, | \, 2) \to G
\end{equation*}
defined in the following way:
\begin{equation} \label{eq:epi-1}
\begin{split}
\varphi(g_1)&=x, \; \; \varphi(g_2)=xy, \; \; \varphi(g_3)=x, \; \;
\psi(g_4)=xy^2, \; \; \psi(g_5)=y, \\
\psi(\ell_1)&=y^2, \; \; \psi(h_1)=y, \; \; \psi(h_2)=x.
\end{split}
\end{equation}
By Riemann Existence Theorem (\cite[Proposition 1.3]{MiPol10}) they
induce two $G$-coverings
\begin{equation*}
f_1 \colon C_1 \to \mathbb{P}^1 \cong C_1/G, \quad f_2 \colon C_2
\to E \cong C_2 /G,
\end{equation*}
where $E$ is an elliptic curve, $g(C_1)=4$ and $g(C_2)=3$.
Moreover, $f_1$ is branched at five points with branching orders
$2$, $2$, $2$, $2$, $4$, whereas $f_2$ is branched at one point
with branching order $2$.

By \eqref{eq:epi-1} it follows that the unique element of $G$,
different from the identity, which acts with fixed points on both
$C_1$ and $C_2$ is $y^2$. Denoting by $|\textrm{Fix}_{C_i}(g)|$ the
number of fixed points of $g \in G$ on the curve $C_i$, we obtain
\begin{equation*}
|\textrm{Fix}_{C_1}(y^2)|=2, \quad |\textrm{Fix}_{C_2}(y^2)|=4,
\end{equation*}
so we have $4 \cdot 2=8$ points in $C_1 \times C_2$ whose
stabilizer is non-trivial (it is isomorphic to $\langle y^2
\rangle$). The $G$-orbit of each of them has cardinality
$|G|/o(y^2)=4$, hence we have exactly two singular points in $T =
(C_1 \times C_2) /G$. More precisely, since $\langle y^2 \rangle$
has order $2$, it follows
\begin{equation*}
\textrm{Sing}\, T = 2 \times \frac{1}{2}(1, \, 1).
\end{equation*}
Let $\lambda \colon X \to T$ be the minimal resolution of
singularities of $T$; then $X$ is a surface of general type, whose
numerical invariants can be computed by using \cite[Proposition
5.1]{MiPol10}; we obtain
\begin{equation*}
p_g(X)=q(X)=1, \quad K_X^2=6.
\end{equation*}
The isotrivial fibration $\pi_2 \colon T \to C_2/G \cong E$ yields,
after composition with $\lambda$,  an isotrivial fibration $\alpha
\colon X \to E$, which is the Albanese morphism of $X$.  The
fibration $\alpha$ has a unique singular fibre $F = 2Y + Z_1 + Z_2$,
where the $Z_i$ are two disjoint $(-2)$ curves and $YZ_i=1$, see
Figure \ref{fig-ex-1}.
\begin{figure}[H]
\begin{center}
\includegraphics*[totalheight=7cm]{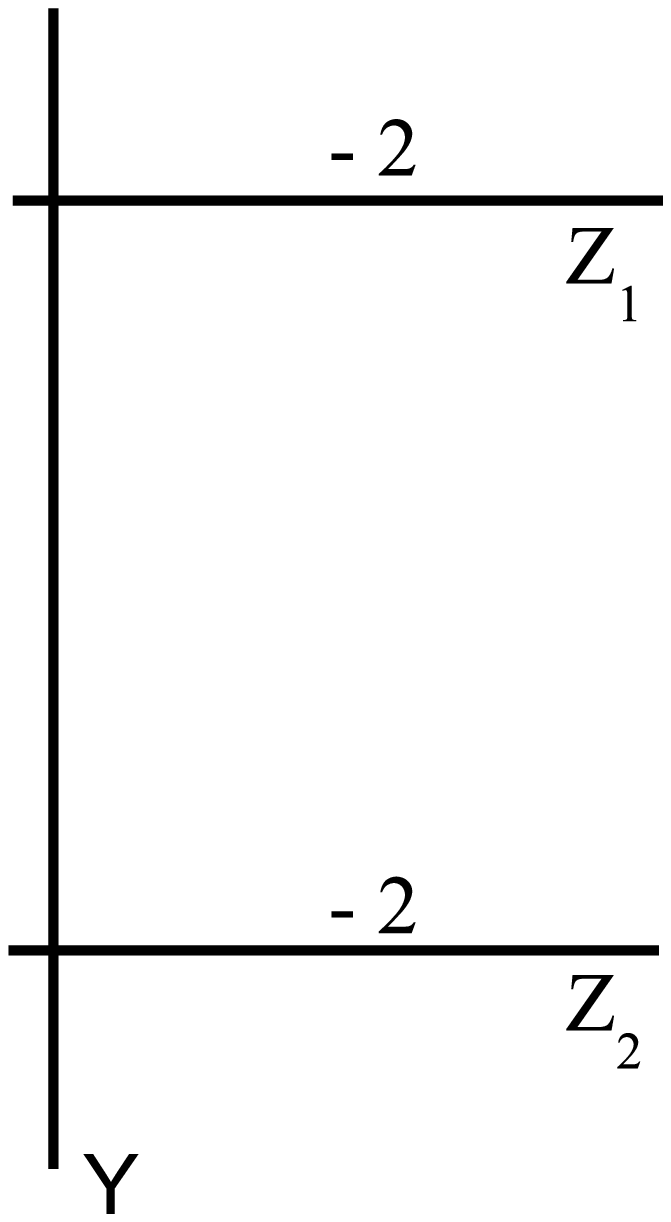}
\caption{The unique singular fibre of $\alpha \colon X \to E$ in
Example \ref{ex.1} and Example \ref{ex.4}} \label{fig-ex-1}
\end{center}
\end{figure}
Since $E$ is an elliptic curve, any $(-1)$-curve of $X$ must be a
component of $F$. On the other hand, we have $K_X F = 2g(C_1)-2=6$
and $F^2=0$, hence $K_XY=3$ and $Y^2=-1$. Thus $Y$ is not a
$(-1)$-curve, so $X$ is a minimal surface of general type satisfying
\begin{equation*}
K_X^2=8 \chi(\mathcal{O}_X)-2
\end{equation*}
and $T$ is the canonical model of $X$. This example shows that
inequality \eqref{eq:main-1} in Theorem \ref{thm:main} is sharp.
Notice that $K_X$ is \emph{not} ample, since $K_XZ_i=0$.
\end{example}

\begin{example} \label{ex.2}
This example can be found in \cite[Section 5]{MiPol10}.

Take as $G$ the semi-direct product $\mathbb{Z}/3 \mathbb{Z} \ltimes
(\mathbb{Z}/4 \mathbb{Z})^2$ whose presentation is
\begin{equation*}
G=\langle x, \, y, \,  z \; | \; x^3=y^4=z^4=1, \, [y, \,z]=1, \,
xyx^{-1}=z, \, xzx^{-1}=(yz)^{-1} \rangle.
\end{equation*}
There are two epimorphisms of groups
\begin{equation*}
\varphi \colon \Gamma(0 \, | \, 3^2, \, 4) \to G, \quad \psi
\colon \Gamma(1 \, | \, 4) \to G
\end{equation*}
defined in the following way:
\begin{equation} \label{eq:epi-2}
\begin{split}
\varphi(g_1)&=x, \; \; \varphi(g_2)=x^2y^3, \; \; \varphi(g_3)=y, \\
\psi(\ell_1)&=y, \; \; \psi(h_1)=x, \; \psi(h_2)=xyxy^2.
\end{split}
\end{equation}
By Riemann Existence Theorem they induce two $G$-coverings
\begin{equation*}
f_1 \colon C_1 \to \mathbb{P}^1 \cong C_1/G, \quad f_2 \colon C_2
\to E \cong C_2 /G,
\end{equation*}
where $E$ is an elliptic curve, $g(C_1)=3$ and $g(C_2)=19$.
Moreover, $f_1$ is branched at three points with branching orders
$3$, $3$, $4$, whereas $f_2$ is branched at one point with
branching order $4$.

By \eqref{eq:epi-2} it follows that the nontrivial elements of $G$
having fixed points on $C_1 \times C_2$ are precisely those in the
set
\begin{equation*}
\Sigma = \bigcup_{\sigma \in G } \langle \sigma y \sigma ^{-1}
\rangle \setminus \{1\}
\end{equation*}
and the elements of order $4$ in $\Sigma$ are $\{y, \, z, \, y^3z^3,
\, y^3, z^3, yz \}$. The product surface $C_1 \times C_2$ contains
exactly $48$ points with nontrivial stabilizer, and for each of them
the $G$-orbit has cardinality $|G|/o(y)=12$, thus $T=(C_1 \times
C_2)/G$ contains $4$ singular points.

Moreover the conjugacy class of $y$ in $G$ is $\{y, \, z, \,
y^3z^3 \}$, whereas the conjugacy class of $y^3$ is $\{y^3, \,
z^3, \, yz \}$. Therefore, for any $h \in \Sigma$ with $o(h)=4$,
one has that $h$ is not conjugate to $h^{-1}$ in $G$. Looking at
the local action of $G$ around each of the fixed points, this
implies that
\begin{equation*}
\textrm{Sing} \, T = 4 \times \frac{1}{4}(1, \, 1).
\end{equation*}
Let $\lambda \colon X \to T$ be the minimal resolution of
singularities of $X$; then $X$ is a surface of general type whose
invariants are
\begin{equation*}
p_g(X)=q(X)=1, \quad K_X^2=2.
\end{equation*}
The isotrivial fibration $\pi_2 \colon T \to C_2/G \cong E$ yields,
after composition with $\lambda$, an isotrivial fibration $\alpha
\colon X \to E$ which is the Albanese morphism of $X$. The
isotrivial fibration $\alpha$ has a unique singular fibre $F = 4Y +
A_1+A_2+A_3+A_4$, where the $A_i$ are disjoint smooth rational
curves such that $A_i^2=-4$ and $YA_i=1$ (see Figure
\ref{fig-ex-2}).
\begin{figure}[H]
\begin{center}
\includegraphics*[totalheight=7 cm]{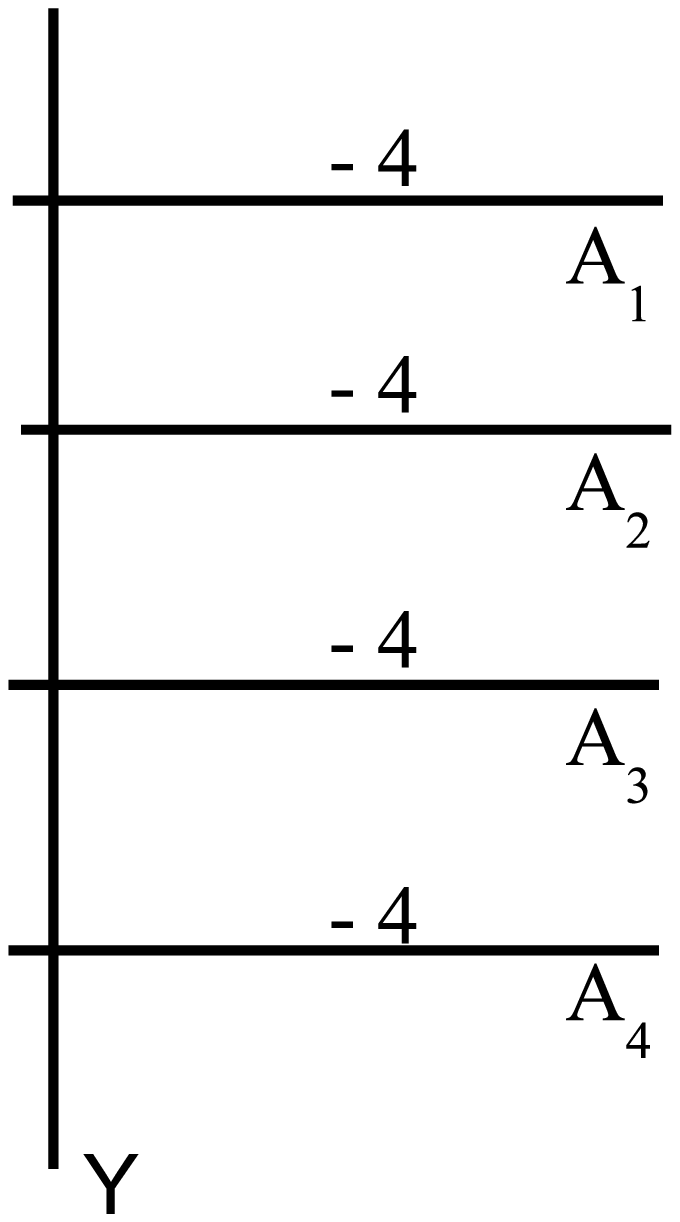}
\caption{The unique singular fibre of $\alpha \colon X \to E$ in
Example \ref{ex.2}} \label{fig-ex-2}
\end{center}
\end{figure}
We have $K_XF=2g(C_1)-2=4$ and $F^2=0$, so we deduce $K_X Y =
Y^2=-1$. Thus $Y$ is a $(-1)$-curve, which is necessarily the
unique $(-1)$-curve in $X$. Since $Y$ is contained in a fibre of
$\alpha \colon X \to E$, after blowing down $Y$ we obtain another
isotrivial fibration $\alpha_m \colon X_m \to E$ such that
\begin{equation*}
p_g(X_m)=q(X_m)=1, \quad K_{X_m}^2=3,
\end{equation*}
that is
\begin{equation} \label{eq:ex2-Xm}
K_{X_m}^2= 8 \chi(\mathcal{O}_{X_m})-5.
\end{equation}
The isotrivial fibration $\alpha_m$ contains a unique singular
fibre, namely the image of $F$ in $X_m$, that is illustrated in
Figure \ref{fig-ex-2-min}.
\begin{figure}[H]
\begin{center}
\includegraphics*[totalheight=7 cm]{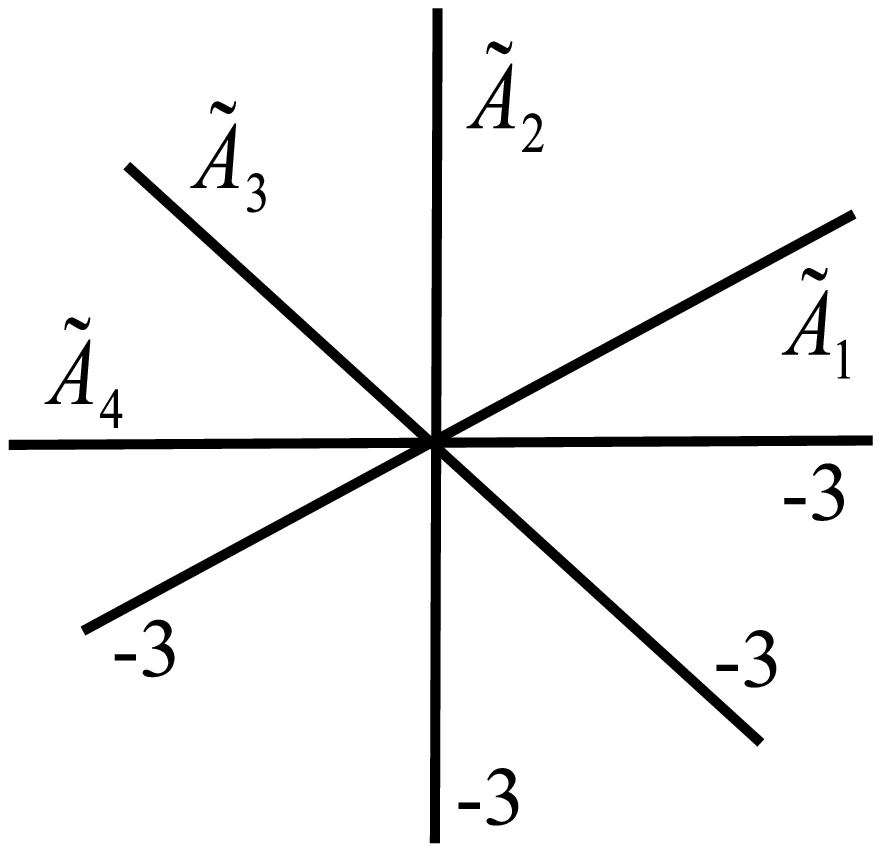}
\caption{The unique singular fibre of $\alpha_m \colon X_m \to E$ in
Example \ref{ex.2}} \label{fig-ex-2-min}
\end{center}
\end{figure}
Here each $\widetilde{A}_i$ is a smooth rational curve with
self-intersection $-3$; in particular $X_m$ contains neither
$(-1)$-curves nor $(-2)$-curves. This shows that $X_m$ is a minimal
model and that $K_{X_m}$ is ample, so \eqref{eq:ex2-Xm} implies that
inequality \eqref{eq:main-2} in Theorem \ref{thm:main} is sharp.
\end{example}

The next two examples show that, if $K_X$ is not ample, then both
cases $K_X^2= 8 \chi (\mathcal{O}_X)-3$ and $K_X^2= 8 \chi
(\mathcal{O}_X)-4$ may actually occur.

\begin{example} \label{ex.3} We exhibit an isotrivial fibration with $K_X$ not
ample and $K_X^2= 8 \chi (\mathcal{O}_X)-3$; this example appears
in \cite[Section 5]{MiPol10}.

Take  $G= D_{12}$, the finite dihedral group of order $12$, with
presentation
\begin{equation*}
G= \langle x, \, y \, | \, x^2=y^6=1, \, xy=y^5x \rangle.
\end{equation*}
There are two epimorphisms of groups
\begin{equation*}
\varphi \colon \Gamma(0 \, | \, 2^3, \, 6) \to G, \quad \psi \colon
\Gamma(1 \, | \, 3) \to G
\end{equation*}
defined in the following way:
\begin{equation} \label{eq:epi-3}
\begin{split}
\varphi(g_1)&=x, \; \; \varphi(g_2)=xy^2, \; \; \varphi(g_3)=y^3, \;
\;
\psi(g_4)=y, \\
\psi(\ell_1)&=y^2, \; \; \psi(h_1)=x, \; \; \psi(h_2)=y.
\end{split}
\end{equation}
By Riemann Existence Theorem they induce two $G$-coverings
\begin{equation*}
f_1 \colon C_1 \to \mathbb{P}^1 \cong C_1/G, \quad f_2 \colon C_2
\to E \cong C_2 /G,
\end{equation*}
where $E$ is an elliptic curve, $g(C_1)=3$ and $g(C_2)=5$. Moreover,
$f_1$ is branched at four points with branching orders $2$, $2$,
$2$, $6$, whereas $f_2$ is branched at one point with branching
order $3$. By \eqref{eq:epi-3} it follows that the nontrivial
elements of $G$ having fixed points on $C_1 \times C_2$ are
precisely those in the set
\begin{equation*}
\Sigma = \bigcup_{\sigma \in G } \langle \sigma y^2 \sigma ^{-1}
\rangle \setminus \{1\}= \{ y^2, \, y^4 \}.
\end{equation*}
The quotient surface $C_1 \times C_2$ contains exactly $8$ points
with nontrivial stabilizer, and for each of them the $G$-orbit has
cardinality $|G|/o(y^2)=4$, thus $T=(C_1 \times C_2)/G$ contains two
singular points. Since $y^2$ is conjugate to $y^4$ in $G$, looking
at the local action of $G$ around each of the fixed points one
obtains
\begin{equation*}
\textrm{Sing} \, T =  \frac{1}{3}(1, \, 1) + \frac{1}{3}(1, \, 2).
\end{equation*}
Let $\lambda \colon X \to T$ be the minimal resolution of
singularities of $X$; then $X$ is a surface of general type whose
invariants are
\begin{equation*}
p_g(X)=q(X)=1, \quad K_X^2=5,
\end{equation*}
that is
\begin{equation*}
K_X^2= 8 \chi (\mathcal{O}_X)-3.
\end{equation*}
The isotrivial fibration $\pi_2 \colon T \to C_2/G \cong E$ yields,
after composition with $\lambda$, an isotrivial fibration $\alpha
\colon X \to E$ which is the Albanese morphism of $X$. The
isotrivial fibration $\alpha$ has a unique singular fibre $F = 3Y +
A+2B_1+B_2$, where $A$ is a $(-3)$-curve and the $B_i$ are
$(-2)$-curves, see Figure \ref{fig-ex-3}.
\begin{figure}[H]
\begin{center}
\includegraphics*[totalheight=7 cm]{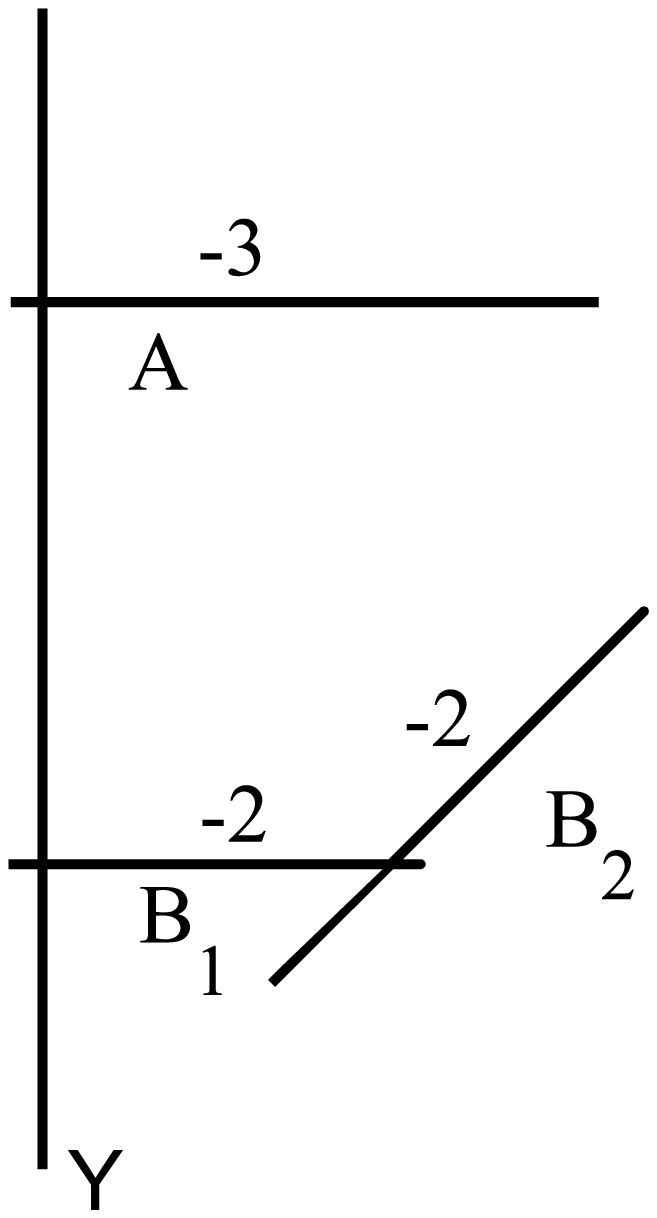}
\caption{The unique singular fibre of $\alpha \colon X \to E$ in
Example \ref{ex.3}} \label{fig-ex-3}
\end{center}
\end{figure}
Using $K_XF=2g(C_1)-2=4$ and $F^2=0$ one obtains $K_XY=1$ and
$Y^2=-1$, hence $Y$ is not a $(-1)$-curve and $X$ is a minimal
model. Notice that $K_X$ is not ample, since $X$ contains the
$(-2)$-curves $B_i$.
\end{example}

\begin{example} \label{ex.4} We exhibit an isotrivial fibration with $K_X$ not
ample and $K_X^2= 8 \chi (\mathcal{O}_X)-4$; this example appears
in \cite[Section 6]{Pol09}.

Take  $G= \mathbb{Z} / 2 \mathbb{Z} \times \mathbb{Z}/2 \mathbb{Z}$,
 with presentation
\begin{equation*}
G= \langle x, \, y \, | \, x^2=y^2=1, \, [x, \, y]=1 \rangle.
\end{equation*}
There are two epimorphisms of groups
\begin{equation*}
\varphi \colon \Gamma(0 \, | \, 2^5) \to G, \quad \psi \colon
\Gamma(1 \, | \, 2^2) \to G
\end{equation*}
defined in the following way:
\begin{equation} \label{eq:epi-4}
\begin{split}
\varphi(g_1)&=x, \; \; \varphi(g_2)=y, \; \; \varphi(g_3)=xy, \; \;
\psi(g_4)=xy, \; \; \varphi(g_5)=xy \\
\psi(\ell_1)&=x, \; \; \psi(\ell_2)=x, \; \;  \psi(h_1)=y, \; \;
\psi(h_2)=y.
\end{split}
\end{equation}
By Riemann Existence Theorem they induce two $G$-coverings
\begin{equation*}
f_1 \colon C_1 \to \mathbb{P}^1 \cong C_1/G, \quad f_2 \colon C_2
\to E \cong C_2 /G,
\end{equation*}
where $E$ is an elliptic curve, $g(C_1)=2$ and $g(C_2)=3$. Moreover,
$f_1$ is branched at five points, $f_2$ is branched at two points
and all the branching orders are equal to $2$. By \eqref{eq:epi-4}
it follows that the unique nontrivial element of $G$ having fixed
points on $C_1 \times C_2$ is $x$. The product surface $C_1 \times
C_2$ contains exactly $8$ points with nontrivial stabilizer, and for
each of them the $G$-orbit has cardinality $|G|/o(x)=2$, thus
$T=(C_1 \times C_2)/G$ contains four singular points. More
precisely, since $\langle x \rangle$ has order $2$,  we have
\begin{equation*}
\textrm{Sing} \, T =  4 \times \frac{1}{2}(1, \, 1).
\end{equation*}
Let $\lambda \colon X \to T$ be the minimal resolution of
singularities of $X$; then $X$ is a surface of general type whose
invariants are
\begin{equation*}
p_g(X)=q(X)=1, \quad K_X^2=4.
\end{equation*}
 The isotrivial fibration $\pi_2 \colon T \to C_2/G \cong E$
yields, after composition with $\lambda$, an isotrivial fibration
$\alpha \colon X \to E$ which is the Albanese morphism of $X$. The
isotrivial fibration $\alpha$ has two singular fibres $F_1$ and
$F_2$, which look like the singular fibre in Figure
 \ref{fig-ex-1}, i.e. they are of the form $F_i=2Y_i+Z_{1i}+Z_{2i}$, where
 the $Z_{ij}$ are disjoint $(-2)$-curves.
 Using $2=K_XF_i=2K_XY_i$ and
$F_i^2=0$ one obtains $K_XY_i=1$ and $Y_i^2=-1$. In particular
$Y_i$ is not a $(-1)$-curve, so $X$ is a minimal surface of
general type satisfying
\begin{equation*}
K_X^2= 8 \chi (\mathcal{O}_X)-4.
\end{equation*}
Notice that $K_X$ is not ample, since $X$ contains the four
$(-2)$-curves $Z_{ij}$.
\end{example}

Looking at Theorem \ref{thm:main}, one might ask whether any
isotrivially fibred surface $X$ with $K_X^2= 8 \chi
(\mathcal{O}_X)-5$ has ample canonical class. The answer is
negative, as shown by the following example that can be found in
\cite[Section 6]{MiPol10}

\begin{example} \label{ex.5}
Take $G=D_{3, 7,2}$, namely the metacyclic group of order $21$
whose presentation is
\begin{equation*}
G=\langle x, \, y \; | \; x^3=y^7=1, \, xyx^{-1}=y^2 \rangle.
\end{equation*}
There are two epimorphisms of groups
\begin{equation*}
\varphi \colon \Gamma(0 \, | \, 3^2, \, 7) \to G, \quad \psi \colon
\Gamma(1 \, | \, 7) \to G
\end{equation*}
defined in the following way:
\begin{equation} \label{eq:epi-5}
\begin{split}
\varphi(g_1)&=x^2, \; \; \varphi(g_2)=xy^6, \; \; \varphi(g_3)=y, \\
\psi(\ell_1)&=y, \; \; \psi(h_1)=y, \; \; \psi(h_2)=x.
\end{split}
\end{equation}
By Riemann Existence Theorem they induce two $G$-coverings
\begin{equation*}
f_1 \colon C_1 \to \mathbb{P}^1 \cong C_1/G, \quad f_2 \colon C_2
\to E \cong C_2 /G,
\end{equation*}
where $E$ is an elliptic curve, $g(C_1)=3$ and $g(C_2)=10$.
Moreover, $f_1$ is branched at three points with branching orders
$3$, $3$, $7$, whereas $f_2$ is branched at one point with branching
order $7$.

By \eqref{eq:epi-5} it follows that the nontrivial elements of $G$
having fixed points on $C_1 \times C_2$ are precisely those in the
set
\begin{equation*}
\Sigma = \bigcup_{\sigma \in G } \langle \sigma y \sigma ^{-1}
\rangle \setminus \{1\}= \{y, \, y^2, \,  y^3, \, y^4, \, y^5, \,
y^6 \}.
\end{equation*}
The product surface $C_1 \times C_2$ contains exactly $9$ points
with non-trivial stabilizer, and for each of them the $G$-orbit has
cardinality $|G|/o(y)=3$; thus $T=(C_1 \times C_2)/G$ contains $3$
singular points.

Moreover, the conjugacy class of $y$ in $G$ is $\{y, \, y^2, \,
y^4\}$ and the conjugacy class of $y^3$ is $\{y^3, \, y^5,
 \, y^6\}$. In particular, every element in $\Sigma$ is conjugate to
 its inverse. Looking at the local action of $G$ around
each of the fixed points, this implies that
\begin{equation*}
\textrm{Sing} \, T = \frac{1}{7}(1, \, 1) + \frac{1}{7}(1, \, 2) +
\frac{1}{7}(1, \, 4).
\end{equation*}
Let $\lambda \colon X \to T$ be the minimal resolution of
singularities of $X$; then $X$ is a surface of general type, whose
invariants are
\begin{equation*}
p_g(X)=q(X)=1, \quad K_X^2=1.
\end{equation*}
The isotrivial fibration $\pi_2 \colon T \to C_2/G \cong E$
yields, after composition with $\lambda$, an isotrivial fibration
$\alpha \colon X \to E$, which is the Albanese morphism of $X$.
The fibration $\alpha$ has a unique singular fibre $F = 7Y +
4A_1+A_2+2B_1+B_2+C$, which looks as in Figure \ref{fig-ex-5}
below.
\begin{figure}[H]
\begin{center}
\includegraphics*[totalheight=7cm]{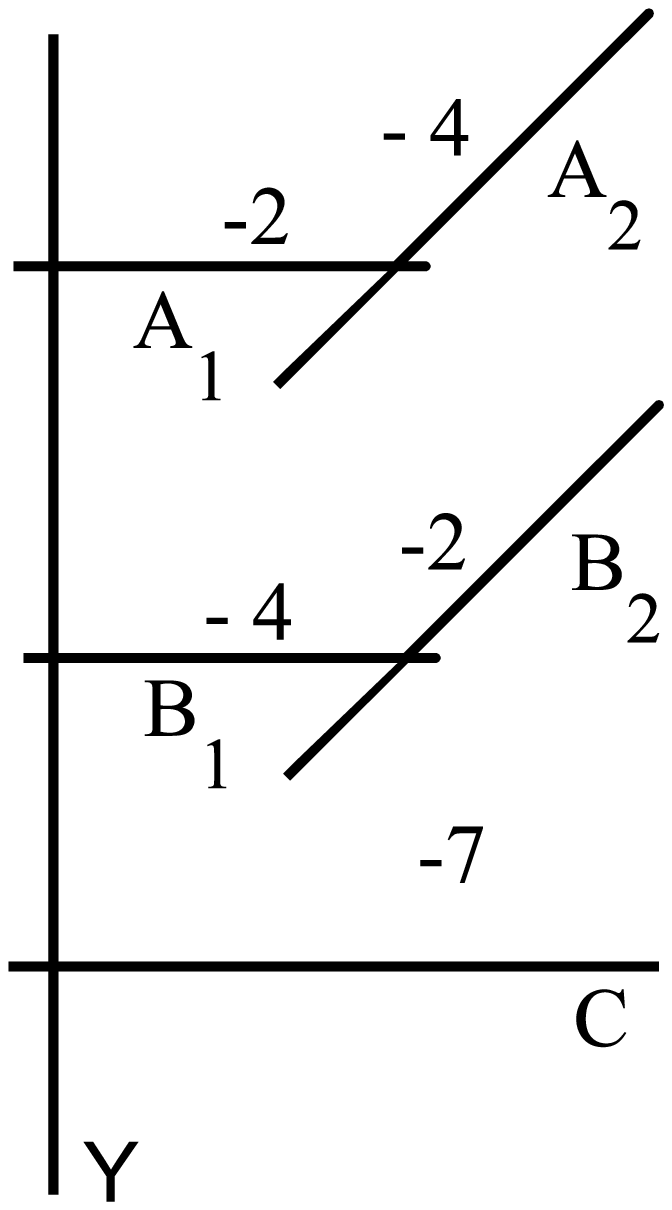}
\caption{The unique singular fibre of $\alpha \colon X \to E$ in
Example \ref{ex.5}} \label{fig-ex-5}
\end{center}
\end{figure}
Here $A_i$, $B_i$ and $C$ are smooth rational curves, and the
integer over each of them denotes as usual the corresponding
self-intersection.

Using $K_XF=4$ and $F^2=0$ we obtain $K_XY=Y^2=-1$, hence $Y$ is the
unique $(-1)$-curve in $X$. The minimal model $X_m$ of $X$ is
obtained by first contracting $Y$ and then the image of $A_1$. Since
$Y$ is contained in a fibre of $\alpha \colon X \to E$, we obtain
another isotrivial fibration $\alpha_m \colon X_m \to E$ such that
\begin{equation} \label{eq:ex3-Xm}
p_g(X_m)=q(X_m)=1, \quad K_{X_m}^2=3,
\end{equation}
that is
\begin{equation*}
K_{X_m}^2= 8 \chi(\mathcal{O}_{X_m})-5.
\end{equation*}
However, the canonical class $K_{X_m}$ is \emph{not} ample, as $X_m$
contains two $(-2)$-curves, namely the images of $B_1$ and $B_2$.
\end{example}

\bigskip
\bigskip

Francesco Polizzi\\
Dipartimento di Matematica e Informatica \\ Universit\`{a}
della Calabria, Cubo 30B  \\ 87036 Arcavacata di Rende, Cosenza (Italy)\\
\emph{E-mail address:} \verb|polizzi@mat.unical.it| \\ \\

\end{document}